\documentclass{article}

\usepackage{geometry}
\usepackage{lmodern} 
\usepackage{wrapfig}
\usepackage{colortbl}
\usepackage{amsmath}
\usepackage{amssymb}
\usepackage{mathrsfs}
\usepackage{amsthm}
\usepackage{makeidx}

\usepackage{multirow}
\usepackage{caption}
\usepackage{arydshln}

\newtheorem{thm}{Theorem}[section]

\usepackage{graphicx}

\usepackage{float}

\floatstyle{ruled}
\newfloat{triangle}{!ht}{tar}
\floatname{triangle}{arithmetic triangle}

\usepackage{ulem}

\usepackage{url}
\date{} 

\title{Distribution function}

\author{Robert Tremblay\\
				Boucherville,
				Canada (Qu\'ebec),\\				
				\texttt{roberttremblay02@videotron.ca}}

\begin{document}


\maketitle

\begin{abstract}
From a known result of diophantine equations of the first degree with $2$ unknowns we simply find the results of the distribution function of the sequences of positive integers generated by the functions at the origin of the $3x+1$ and $5x+1$ problems. 
\end{abstract}

\section{Introduction}
 
We demonstrate and analyze the properties of the distribution function $F(k)$ introduced by Riho Terras~\cite{terras} in $1976$ and taken up by several authors including Lagarias~\cite{lagarias} using a whole new approach. Terras proves that this function is well defined and has very interesting properties without ever, however, solving the conjecture linked to the 3x+1 problem. Nevertheless, it remains probably one of the biggest advances in this quest. In our opinion, the path used by Terras is the most complete. However, this requires 5 definitions, 11 theorems, 8 corollaries and 4 propositions. In our paper we only use three theorems to reach the same results.      

First, we rediscover the periodicity theorem on the distribution of the trajectories generated by the iteration of the function at the origin of the $3x+1$ problem, which also applies to the $5x+1$ problem. Terras qualified this result as a remarkable phenomenon of periodicity. Subsequently, we find exactly the results presented by Terras.

Detailed analysis of the function and its various properties will lead us to a somewhat unexpected conclusion about the conjecture.

\section{Functions $T_3$ and $T_5$}

Mappings can be define on integers represented by functions such that each element of the set $\mathbb{Z}$ is connected to a single element of this set. The iterative application of these functions produces a sequence of integers called trajectories. 

Let

\begin{equation*}
(n, f(n),f^{(2)}(n),f^{(3)}(n), \cdots, f^{(i)}(n), \cdots) ,
\end{equation*}

 with $f^{(i+1)}(n) = f\{f^{(i)}(n)\}$, $i =0, 1, 2, 3, \cdots$ and $f^{(0)}(n) = n$, a trajectory generates by a function f on an integer $n$.

A sequence of integers forms a loop when there exists a number of iterations $k \ge 1$ such that 

\begin{equation}
f^{(k)}(n) = n. 
\label{condition_cycle}
\end{equation}

If all integers in the sequence are different two by two, we have by definition a cycle of length $p = k$, so the trajectory $(n, f(n),f^{(2)}(n),f^{(3)}(n), \cdots, f^{(k-1)}(n))$. Generally, we note the trajectory characterizing a cycle starting with the smallest integer. 

There are a multitude of functions that have these properties. The function $g(n)$ giving rise to the original Collatz problem and the $3x+1$ function $T_3(n)$~\cite{lagarias}, the $5x+1$ function $T_5(n)$ and the accelerated $3x+1$ function~\cite{kontorovich_lagarias}, are some examples. Except for the last function, the others come from a group called Generalized $3x+1$ Mappings~\cite{matthews}.

The two functions dealt with in this paper are defined by 

\begin{equation}
T_{m_i}(n) = \left\lbrace
\begin{array}{ll}
\frac{n}{2} & \mbox{, if $n\equiv0\pmod{2}$},\\
\\
\frac{m_in+1}{2} & \mbox{, if $n\equiv1\pmod{2}$}\\
\end{array}
\right.
\label{mx+1_function}
\end{equation}

with $m_3 = 3$ for the $3x+1$ problem and , $m_5 = 5$ for the $5x+1$ problem.

The general expression giving the result of $k$ iterations of the function $T_{m_i}$, which we will simply call $T$, on an integer $n$ is 

\begin{equation}
	T^{(k)}(n) = \lambda_{k_1,k_2}n + \rho_k(n),  
\label{g_k}
\end{equation}

where

\begin{equation}
	\lambda_{k_1,k_2} = \left(\frac{1}{2}\right)^{k_1}\left(\frac{m_i}{2}\right)^{k_2} 
\label{lambda_k1_k2}
\end{equation}

and

\begin{equation}
	k = k_1 + k_2,
\label{k1_plus_k2}
\end{equation}

with $k_1$ the number of transformations of the form $n/2$ and $k_2$, transformations of the form $(m_in + 1)/2$.

Unlike parameter $\lambda_{k_1,k_2}$, $\rho_k(n)$ depend on the order of application of the transformations.

Let $n$ and $T^{(k)}$ be replaced by the variables $x$ and $y$,

\begin{equation*}
	2^k\rho_k(n) = 2^ky - 1^{k_1}m_i^{k_2}x.  
\label{eq_diophantine}
\end{equation*}

In this form we have a diophantine equation of first degree at two unknowns, 

\begin{equation}
	c = by - ax,  
\label{eq_diophantine_general}
\end{equation}

where 

\begin{equation}
	a = 1^{k_1}m_i^{k_2}, \phantom{12} b = 2^k \phantom{12} and \phantom{12} c = 2^k\rho_k(n).  
\label{par_a_b_c}
\end{equation}

Depending on the new parameters $a$ and $b$ the parameter $\lambda$ (equation \ref{lambda_k1_k2}) becomes 

\begin{equation}
	\lambda_{a,b} = \left(\frac{a}{b}\right). 
\label{lambda_a_b}
\end{equation}

From a well-known result of diophantine equations theory we have the theorem

\begin{thm}
Let the diophantine equation $c = by - ax$ of first degree at two unknowns. If the coefficients $a$ and $b$ of $x$ and $y$ are prime to one another (if they have no divisor other than $1$ and $-1$ in common), this equation admits a infinity of solutions to integer values. If $(x_0, y_0)$ is a specific solution, the general solution will be $(x = x_0 + bq, y = y_0 + aq)$, where $q$ is any integer, positive, negative or zero.
	\label{equationDiophantine}
\end{thm}

\textsl{Proof}

References : Bordell\`es~\cite{bordelles}. $\blacksquare$

We may to assign to every integer of a trajectory generates by the function $T(n)$ a number $t_j = 0$ if $T^{(j)}(n)$ is even, and $t_j = 1$ if it is odd. Then, the iterative application of the function $T$ to an integer $n$ give a diadic sequence $w_l$ of $1$ and $0$  


\begin{equation*}
w_l = (t_0, t_1, t_2, t_3, \cdots, t_j, \cdots, t_{l-1}), \phantom{1,2} with \phantom{1,2} l \ge 1.
\end{equation*} 

Diadic sequences correspond to what Lagarias called the parity vectors~\cite{lagarias}.

For a given length $l$ there are $2^l$ different diadic sequences $w_l$ of $0$ and $1$. 

The representation of the trajectories in terms of $t_{j}$ leads to an important theorem which makes it possible to bring out an intrinsic property, namely the \textsl{periodicity}. This property has already been observed by Terras~\cite{terras} and Everett~\cite{everett} concerning the process of iterations of the function $T_3(n)$ generating the problem $3x+1$, and appears in a theorem which they have demonstrated by induction. We will prove it differently, using the previous theorem.


\begin{thm}
	All diadic sequences $w_l$ of length $l = k \ge 1$ generated by any $2^l$ consecutive integers are different and are repeated periodically.
	\label{periodicity}
\end{thm}

\textsl{Proof}

Let $k = l \ge 1$ the number of iterations applied to a given integer $n$. The trajectories 

\begin{flushleft}
	$\phantom{1,2,3,4} (T^{(0)}(n),T^{(1)}(n))$ 
	
	$\phantom{1,2,3,4} (T^{(0)}(n),T^{(1)}(n),T^{(2)}(n))$ 
	
	$\phantom{1,2,3,4} \cdots$ 
	
	$\phantom{1,2,3,4} (T^{(0)}(n),T^{(1)}(n)), \cdots, T^{(k)}(n))$
\label{trajectories}
\end{flushleft}

correspond respectively to the diadic sequences 

\begin{flushleft}
	$\phantom{1,2,3,4} w_1 = (t_0)$ 
	
	$\phantom{1,2,3,4} w_2 = (t_0, t_1)$ 
	
	$\phantom{1,2,3,4} \cdots$
	
	$\phantom{1,2,3,4} w_{l=k} = (t_0, t_1, \cdots, t_{k-1})$ .
\label{sequences}
\end{flushleft}

For a given number $l$ we have $2^l$ different diadic sequences $w_l$ possible. 

According to theorem~\ref{equationDiophantine}, each of the $2^l$ diadic sequences will be performed for $k = l$. Indeed, the $0$ and the $1$ of these sequences correspond to the operations on the even and odd integers. We build $2^k$ different diophantine equations characterized by $2^k$ different combinations of the parameters $a$, $b$ and $c$, whose solutions will be given by $(x = x_0 + 2^kq, y = y_0 + m_i^{k_2}q)$. Therefore, all the integers $x_0 + 2^kq$ starting a trajectory of length $k + 1$ correspond to the same sequence $w_k$. In a sequence of $2^k$ consecutive integers, each integer must start a different sequence $w_k$, otherwise the $2^k$ different diadic sequences will not be performed. $\blacksquare$


We will use another property of the diophantine equations generated by functions like $T_3$ and $T_5$.

\begin{thm}
Let the trajectories of the integers (of length L) that are connected to each other by the operations $n/2$ or $(m_in + 1)/2$. The diophantine equation connecting the first integer $x$ and the last integer $y$ of a sequence can be expressed in the general form $c = by - ax$ where the parameters $a$, $b$ and $c$, always positive, depend on the operations themselves and in which orders they are applied. If $b < a$($\lambda > 1$), $x < y$ and, if $b > a$($\lambda < 1$), $x \ge y$ or $x<y$ when $x$ and $y$ are positive.
	\label{Distribution_x_vs_y}
\end{thm}

\textsl{Proof}

Let $k_1, k_2 = 0, 1, 2, \ldots$ and $k = k_1 + k_2 = L - 1$, with $L \ge 2$.

Then, $a = m_i^{k_2}$, $b = 2^k$ and $c \ge 0$.

As the factors $a$ and $b$ of $x$ and $y$ are prime to one another, the diophantine equation admits a infinity of solutions to integer values. If $(x_0, y_0)$ is a specific solution, the general solution will be $(x = x_0 + bq, y = y_0 + aq)$, where $q$ is any integer, positive, negative or zero.

Let the equation~\ref{g_k} in the form $y = \lambda x + \rho$, with $\lambda = a/b$ and $\rho$ always positive. A quick examination of this equation allows us to state that if $\lambda > 1$ ($b < a$), $x < y$, and if $\lambda < 1$ ($b > a$), $x \ge y$ or $x < y$. Therefore, two cases are possible, so $b < a$ or $b > a$. Now, let us examine these two cases from the diophantine equation and its solutions. 


\underline{First case : $b < a$}
 
Of the diophantine equation $c = by - ax$, as $c$ is always positive and $b < a$, $x$ must always be smaller than $y$ ($x < y$).

\underline{Second case : $b > a$}

Let the general solution

\begin{equation*}
	y = y_0 + aq \phantom{1234} and \phantom{1234} x = x_0 + bq = x_0 + 2^kq,	
\end{equation*}


There are combinations $k = k_1 + k_2$ giving parameters $b$ and $a$ such that $b > a $ ($\lambda < 1$), and maybe integers $x_0 < y_0$ included in the interval $1$ to $2^k$. As $b > a$ and $x_0 < y_0$, beyond a certain value of $q = q_{x > y}$, we will have $x > y$. All integers $x = x_0 + 2^kq$ with $q > q_{x > y}$ will be greater than $y$.

If $x = y$ then $c = x(b - a)$. Since $c$ must always be positive, then $b > a$. $\blacksquare$




For example, for $L = 2$, we have $k = k_1 + k_2 = L - 1 = 1$. Two cases are possible, $k_1 = 1$ and $k_2 = 0$ or, $k_1 = 0$ and $k_2 = 1$. Then, if $m_i = m_3 = 3$, we have $a = 3^{k_2} = 3^0 = 1$ or $a = 3^{k_2} = 3^1 = 3$ and, $b= 2^k = 2$. We write the diophantine equations

\begin{equation}
	0 = 2y - x \phantom{1234} or \phantom{1234} 1 = 2y - 3x.	
\label{ED_k_eq_2}
\end{equation}

where $(b = 2, a = 1)$ ($b > a$) in the first case and $(b = 2, a = 3)$ ($b < a$) in the other case.

The general solutions $(x, y)$ are respectively $(2 + 2q, 1 + 1q)$ with ($x > y$), and $(1  + 2q, 2 + 3q)$ with ($x < y$).

          
\section{Distribution function $F(k)$}

Let us define the distribution function $F(k)$ as

\begin{equation}
	F(k) = \lim\limits_{m \rightarrow \infty} \left( 1/m \right) \mu \{ n \le m \phantom{1} | \phantom{1} \chi(n) \ge k \},
\label{Function_distribution}	
\end{equation} 

where $\mu$ is the number of positive integers $n \le m$ with $m$ that tends towards infinity. $\chi(n)$ is called the "stopping time", and corresponds to the smallest positive integer such that the iterative application ($k$ times) of function $T_3$ (equation~\ref{mx+1_function}) on a integer $n$ gives the result $T_3^{(k)}n<n$.

We can state the $3x+1$ conjecture~\cite{terras, lagarias} as follows:

$3x+1$ CONJECTURE. Every integer $n \ge 2$ has a finite stopping time.

Terras~\cite{terras} proves that the distribution function $F(k)$ is well defined for any value of $k$ and that it tends towards $0$ for $k$ tending towards infinity.

Lagarias~\cite{lagarias} redoes the demonstration using the function we will call $G(k)$,

\begin{equation}
	G(k) = \lim\limits_{x \rightarrow \infty} \frac{1}{x} \phantom{1} \# \{ n : n \le x \phantom{1} and \phantom{1} \sigma(n) \le k \},
\label{Function_distribution_Lagarias}	
\end{equation} 

where $\sigma(n)$ is the "stopping time". This function $G(k)$ is in away almost the reciprocal of the function $F(k)$, and tends towards $1$ when $k$ tends towards infinity. The properties inherent in these functions will be clarified in the following examples.

The application of theorem~\ref{periodicity} on periodicity can be interpreted as follows.

Let $k$ be a number of iterations applied to any $2^k$ consecutive integers. We will have all possible combinations $2^k$ of operations $n/2$ on the even integers and $(3n+1)/2$ on the odd integers of the diadic sequences generated by the function $T_3(n)$ and each combination appears only once. For a given $k$, all the integers $m$ of the form $m = n + 2^kq$ will have the same combination of operations. The distribution of different combinations is then binomial versus the operations.

For example, let $k=1$ and the $2^k=2^1=2$ consecutive positive integers $1$ and $2$. The trajectories of length $k + 1 = 2$ generated by the function $T_3(n)$ will be

\begin{equation*}
	(1,2) \phantom{1} (2,1) \phantom{1} (3,5) \phantom{1} (4,2) \phantom{1} (5,8) \phantom{1} (6,3)  \phantom{1} \cdots,
\label{sequences_2_numbers}	
\end{equation*} 

where we have added the numbers $3$, $4$, $5$ and $6$ after the two consecutive numbers $1$ and $2$ starting the trajectories, so as to bring out the periodicity.

If we use the diadic sequences of the $0$ and $1$ representing respectively the even and odd operations, we will have

\begin{equation*}
	(1) \phantom{1} (0) \phantom{1} (1) \phantom{1} (0) \phantom{1} (1) \phantom{1} (0) \phantom{1} \cdots,
\label{sequences_2_numbers_diadic}	
\end{equation*} 

all repeating periodically for every two consecutive trajectories. This result follows from the fact that the all integers alternate between the even and odd integers. 

We have already writed the diophantine equations for $k = 1$ (equations~\ref{ED_k_eq_2}) which give the first integer $x$ of the trajectory versus the last integer (here the second).

In the first case we have all the trajectories starting with an even positive integer $x$ and ending with a smaller integer $y$ after $1$ iteration. The stopping time is equal to the number of iterations $k = 1$, so $\chi(n = even) = k = 1$. In the second case we have all trajectories starting with an odd positive integer $x$ and ending with a greater integer $y$ after $1$ iteration and, $\chi(n = odd) > k = 1$. The stopping time meets the condition $\chi(n) \ge k = 1$ in two cases and all integers contribute to the distribution function $F(k)$, so $F(k = 1) =  1$. Unlike Terras, we will not count the integers with $\chi = k$ because in these cases, we have reached the condition $T_3^{(k)}n < n$. It will create a slight gap with the results of Terras. Then, the distribution function $F(k)$  with $\chi > k$ instead $\chi \ge k$ really becomes the reciprocal of the function $G(k)$ defined by Lagarias, and the new function $F_{new}(k = 1) = 1/2$. We write

\begin{equation}
	F_{new}(k) = \lim\limits_{m \rightarrow \infty} \left( 1/m \right) \mu \{ n \le m \phantom{1} | \phantom{1} \chi(n) > k \}.
\label{Function_distribution_new}	
\end{equation}



Let another example. Take $k=2$ and the $2^k=2^2=4$ consecutive positive integers $3,4,5$ and $6$. The trajectories of length $k + 1 = 3$ generated by the function $T(n)$ will be

\begin{equation*}
	(3,5,8) \phantom{1} (4,2,1) \phantom{1} (5,8,4) \phantom{1} (6,3,5) \phantom{1} (7,11,17) \phantom{1} (8,4,2)  \phantom{1} \cdots,
\label{sequences_4_numbers}	
\end{equation*} 

where we have added the numbers $7$ and $8$ after the four consecutive numbers $3$, $4$, $5$ and $6$ starting the trajectories, so as to bring out the periodicity. 

The diadic sequences are

\begin{equation*}
	(1,1) \phantom{1} (0,0) \phantom{1} (1,0) \phantom{1} (0,1) \phantom{1} (1,1) \phantom{1} (0,0) \phantom{1} \cdots,
\label{sequences_4_numbers_diadic}	
\end{equation*} 

all repeating periodically for every four consecutive integers starting a trajectory. 

We can write the $2^k = 2^2 = 4$ diophantine equations in the same way as before. But, we will do it differently here. In fact the diadic sequences and the theorem~\ref{Distribution_x_vs_y} we will help to deduce whether or not the stopping time is equal, greater or less than the number of iterations $k = 2$.

In the general case, the parameter $b = 2^k$ and the parameter $a = 3^{k_2} \cdot 1^{k_1} = 3^{k_2}$ with $k$ the total number of iterations, $k_1$ the number of operations on the even integers, and $k_2 = k - k_1$ the number of operations on the odd integers. 

As the first two diadic sequences (table~\ref{stoppingTime_K2}) correspond to the trajectories starting with an even integer, we do not count them in $F(k)$. The third diadic sequence, so $(1,0)$ which is generated by the integers $5 + 4q$, is such that $\chi(5 + 4q) = k = 2$. The fourth diadic sequence, so $(1,1)$ which is generated by the integers $3 + 4q$, is such that $\chi(3 + 4q) > k = 2$. Then, the distribution function $F_{new}(k)$  with $\chi > k$ instead $\chi \ge k$ really becomes $F_{new}(k=2) = 1/4$. The original function $F(k = 2)$ would correspond to $1/2$.
 
And so on for different values of the number of iterations $k$.


\begin{table}[H]
\begin{center}
\begin{tabular}{|c|c|c|c|c|c|c|c|}
	
	\hline
	$diadic$ & $k_1$ & $k_2$ & $b = 2^k$ & $a = 3^{k_2}$ & $b \phantom{1} vs \phantom{1} a$ & $x \phantom{1} vs \phantom{1}  y$ & $stopping$ \\
	
	$sequences$ & & & & & & & $time \phantom{1} \chi(n)$ \\
	\hline
	$(0,0)$ & 2 & 0 & 4 & 1 & $b > a$ & $x > y$ & $-$ \\
	
	$(0,1)$ & 1 & 1 & 4 & 3 & $b > a$ & $x > y$ & $-$ \\
	
	$(1,0)$ & 1 & 1 & 4 & 3 & $b > a$ & $x > y$ & $\chi = k$ \\
	
	$(1,1)$ & 0 & 2 & 4 & 9 & $b < a$ & $x < y$ & $\chi > k$ \\
	\hline
\end{tabular}
\end{center}
\caption{Stopping time for k = 2}
\label{stoppingTime_K2}

\end{table}

For a given $k$, the total number of different trajectories is $b = 2^k$. For a given $0 \le k_2 \le k$ the number of different trajectories is calculated by the binomial coefficients $\left( \begin{array}{c} k \\ k_2 \end{array}\right)$. Binomial coefficients (BC) can be represented in a Pascal triangle (table~\ref{PascalTriangle_BC}),


\begin{table}[H]
\begin{center}
\begin{tabular}{c|cccccccccccc}
	
	$k_2 \setminus k$  & 0 & 1 & 2 & 3 & 4 & 5  & 6  & 7  & 8  & 9 & 10 &$\cdots$ \\
		\hline
	0 & 1 & 1 & 1 & 1 & 1  & 1  & 1  & 1  & 1  & 1 & 1 & $\cdots$ \\
	1 &   & 1 & 2 & 3 & 4  & 5  & 6  & 7  & 8  & 9 & 10 & $\cdots$ \\
	2 &	  &		& 1 & 3 & 6  & 10 & 15 & 21 & 28 & 36 & 45 & $\cdots$ \\
	3	&		&		&   & 1 & 4  & 10 & 20 & 35 & 56 & 84 & 120 & $\cdots$ \\
	4	&		&		&		&   & 1  & 5  & 15 & 35 & 70 & 126 & 210 & $\cdots$ \\
	5	&		&		&		&   &    & 1  & 6  & 21 & 56 & 126 & 252 & $\cdots$ \\
	6	&		&		&		&   &    &    & 1  & 7  & 28 & 84 & 210 & $\cdots$ \\			
	7	&		&		&		&   &    &    &    & 1  & 8  & 36 & 120 & $\cdots$ \\
	8	&		&		&		&   &    &    &    &    & 1  & 9 & 45 & $\cdots$ \\
	9	&		&		&		&   &    &    &    &    &    & 1 & 10 & $\cdots$ \\
	10	&		&		&		&   &    &    &    &    &    &  & 1 & $\cdots$ \\
	$\cdots$	&		&		&		&   &    &    &    &    &    &  &  &  \\
	total	&	1	&	2	&	4	& 8  &  16  &  32  &  64  &  128  &  256  & 512 & 1 024 & $\cdots$ \\
	

\end{tabular}
\end{center}
\caption{Pascal triangle - Binomial coefficients}
\label{PascalTriangle_BC}

\end{table}

If we number each of the rows $i = k_2$ and each of the columns $j =k$, we can write $n(i,j) = BC$, so the different elements of the table. By the properties of the binomial coefficients we have $n(0,k) = 1$ (the top line) and $n(k,k) = 1$ (the bottom diagonal). In addition, all other elements are the result of sum 

\begin{equation}
	n(0 < i < j, j > 1) = n(i-1,j-1) + n(i,j-1).	
\label{somme}
\end{equation}

We use a similar table (table~\ref{PascalTriangle_nbr}) which contain the number of integers $n(i = k_2,j = k)$ by $2^k$ consecutive integers which satisfy the condition that the the stopping time $\chi$ is greater than the number of iterations $k$. We have 


\begin{table}[H]
\begin{center}
\begin{tabular}{c|cccccccccccc}
	
	$k_2 \setminus k$  & 0 & 1 & 2 & 3 & 4 & 5  & 6  & 7  & 8 & 9 & 10 & $\cdots$ \\
		\hline
	0   & 1 & 1\cellcolor[gray]{0.8} & 0 & 0 & 0  & 0  & 0  & 0  & 0  & 0  & 0  & $\cdots$ \\
	1   &   & 1 & 1\cellcolor[gray]{0.8} & 0 & 0  & 0  & 0  & 0  & 0  & 0  & 0  & $\cdots$ \\
	2   &	  &		& 1 & 1 & 1\cellcolor[gray]{0.8}  & 0  & 0  & 0  & 0  & 0  & 0  & $\cdots$ \\
	3	  &		&		&   & 1 & 2  & 2\cellcolor[gray]{0.8}  & 0  & 0  & 0  & 0  & 0  & $\cdots$ \\
	4	  &		&		&		&   & 1  & 3  & 3  & 3\cellcolor[gray]{0.8}  & 0  & 0  & 0  & $\cdots$ \\
	5	  &		&		&		&   &    & 1  & 4  & 7  & 7\cellcolor[gray]{0.8}  & 0  & 0  & $\cdots$ \\
	6	  &		&		&		&   &    &    & 1  & 5  & 12 & 12 & 12\cellcolor[gray]{0.8}  & $\cdots$ \\			
	7	  &		&		&		&   &    &    &    & 1  & 6  & 18 & 30 & $\cdots$ \\
	8	  &		&		&		&   &    &    &    &    & 1  & 7  & 25 & $\cdots$ \\ 
	9	  &		&		&		&   &    &    &    &    &    & 1  & 8  & $\cdots$ \\ 
	10	&		&		&		&   &    &    &    &    &    &    & 1  & $\cdots$ \\ 		
	$\cdots$ 	&		&		&		&   &    &    &    &    &    &    &    & $\cdots$ \\
	total	&	1	&	1	&	1	& 2  &  3  &  4  &  8  &  13  &  19  & 38 & 64 & $\cdots$ \\

\end{tabular}
\end{center}
\caption{Pascal triangle - Number of integers $n(i = k_2,j = k)$ by $2^k$ consecutive integers with $\chi > k$ ($a > b$ and $y > x$)}
\label{PascalTriangle_nbr}

\end{table}


The values in the shaded areas for a given $k$, correspond to the condition $a_k < b_k$ ($\lambda_k < 1$) when $a_{k-1} > b_{k-1}$ ($\lambda_{k-1} > 1$) so, the condition used by Terras and Lagarias to define an admissible vector. Unlike Terras, we do not count these values in the calculation of the total. 


The index $j$ for the columns of the table is the exponent $k$ (the number of iterations) of $2$ in the parameter $b = 2^k$. The index $i$ for the rows is the exponent $k_2$ of $3$ in the parameter $a = 3^{k_2}$. As $k_2$ correspond to the number of operations on the odd integers, this value is in fact the number of $1$ in the diadic sequences and varies of $0$ to $k$. The various data in this table are calculated recursively.

The first data is trivial and indicates that all the integers satisfy the condition $\chi > k$ and this, because the number of iterations is $k = 0$. The case $k = 1$ has ready be analyzed and we have $n(0,1) = 0$ and $n(1,1) = 1$. After $1$ iteration, all positive even integers go to a smaller integer $(n = 0)$ and, all positive odd integers go to a greater integer $(n = 1)$.   

From $k = 2$ we proceed recursively in the calculation of $n(i,k)$. 

We use the principle that each sequence is generated so that the new parameter $b$ (for $k$) is the precedent (for $k - 1$) time $2$, and the new parameter $a$ (for $k_2$) is the precedent (for $k_2 -1$) time $1$ or $3$. 

For example, for $k = 2$, we have two $n$ which precede (for $k = 1$), so $n(0,1) = 0$ and $n(1,1) = 1$. As $n(0,1) = 0$, the sequences starting with a even positive integer for $k = 2$ will not contribute to $F(k)$ and $n(0,2) = 0$. On the other hand, the sequences generated by the integers with $n(1,1) = 1$ can contribute to $n(1,2)$ and $n(2,2)$. The new parameter $b$ will be $b = 2 \cdot 2$ and the new parameter $a$ will be $a = 3 \cdot 1$ or $a = 3 \cdot 3$ (table~\ref{stoppingTime_K2}). In the first case, $b > a$, $x > y$ and $\chi = k$. Then $n(1,2) = 0$. In the second case, $b < a$, $x < y$ and $\chi > k$. Then $n(2,2) = 1$. And so on for different values of $k$. 

We put zeros for $n(i,j)$ when $b > a$.

The sum on the index $i$ of $n(i,k) / 2^k$ for a given $k$ gives the value of the distribution function $F_{new}(k)$ for this number of iterations $k$. Knowing that non-zero values must satisfy inequality $b < a$ ($x < y$), with $b = 2^k$ and $a = 3^{k_2}$, the sum begins with $i = k_2 > k \theta$,  

\begin{equation}
	F_{new}(k) = \sum_{i = 0}^{k} \frac{n(i,k)}{2^k} = \sum_{i > k\theta}^{k} \frac{n(i,k)}{2^k}, \phantom{1} with \phantom{1} \theta = \frac{\ln{2}}{\ln{3}} \simeq 0.63093.
\label{Function_distribution_new_sum}		
\end{equation} 

It is then easy to build the computer programs starting from the recursive function worked out by Terras and by the previous process which makes it possible to fill the table~\ref{PascalTriangle_nbr}. The results of these two programs are compiled in the table~\ref{Distribution function F(k)_3(k)}. 

\begin{table}
\begin{center}


\begin{tabular}{c|c|c|c|c|c|c}
	\hline
	 $k$ &   Terras & new & & $k$ & Terras & new \\
	\hline
		10 &  $7.4219 \times 10^{-2}$ &  $6.25 \times 10^{-2}$ &  & 100 & $2.6396 \times 10^{-4}$ & $2.3868 \times 10^{-4}$\\

		20 &  $2.8591 \times 10^{-2}$ &  $2.6062 \times 10^{-2}$ &  & 200 & $3.3187 \times 10^{-6}$ & $3.0604 \times 10^{-6}$ \\	

	  30 &  $1.1894 \times 10^{-2}$ &  $1.1894 \times 10^{-2}$ &  & 300 & $5.7714 \times 10^{-8}$ & $5.4667 \times 10^{-8}$ \\	
		
		40 &  $6.5693 \times 10^{-3}$ &  $5.8233 \times 10^{-3}$ &  & 400 & $1.2191 \times 10^{-9}$ & $1.1587 \times 10^{-9}$ \\
			
		50 &  $3.5373 \times 10^{-3}$ &  $3.3167 \times 10^{-3}$ &  & 500 & $2.7866 \times 10^{-11}$ & $2.6584 \times 10^{-11}$ \\
		
		60 &  $1.9222 \times 10^{-3}$ &  $1.9222 \times 10^{-3}$ &  & 600 & $6.7168 \times 10^{-13}$ & $6.4455 \times 10^{-13}$ \\
		
		70 &  $1.1644 \times 10^{-3}$ &  $1.0516 \times 10^{-3}$ &  & 700 & $1.5719 \times 10^{-14}$ & $1.5719 \times 10^{-14}$ \\
		
		80 &  $7.0744 \times 10^{-4}$ &  $6.6440 \times 10^{-4}$ &  & 800 & $4.0963 \times 10^{-16}$ & $4.0963 \times 10^{-16}$ \\
		
		90 &  $4.1078 \times 10^{-4}$ &  $4.1078 \times 10^{-4}$ &  & 900 & $1.0837 \times 10^{-17}$ & $1.0837 \times 10^{-17}$ \\
	
	\hline
	\end{tabular}
	
	\end{center}
	\caption{Distribution function $F_3(k)$}
	\label{Distribution function F(k)_3(k)}
	\end{table}

We have also extended the programs to the distribution function $F_5(k)$ generated by the $5x + 1 $ function $T_5$ (table~\ref{Distribution function F(k)_5(k)}).
	
\begin{table}
\begin{center}


\begin{tabular}{c|c|c|c|c|c|c}
	\hline
	 $k$ &   Terras & new & & $k$ & Terras & new \\
	\hline
		10 &  0.2734375 &  0.25976563 &  & 100 & 0.18087772 & 0.18060217\\

		20 &  0.22122192 &  0.22122192 &  & 200 & 0.17688689 & 0.17685114 \\	

	  30 &  0.20572651 &  0.20572651 &  & 300 & 0.17622449 & 0.17621811 \\	
		
		40 &  0.19784735 &  0.19625785 &  & 400 & 0.17607927 & 0.17607775 \\
			
		50 &  0.19116563 &  0.19116563 &  & 500 & 0.17604079 & 0.17604048 \\
		
		60 &  0.18811449 &  0.18811449 &  & 600 & 0.17603033 & 0.17603024 \\
		
		70 &  0.18573498 &  0.18513014 &  & 700 & 0.17602715 & 0.17602715 \\
		
		80 &  0.18317774 &  0.18317774 &  & 800 & 0.17602622 & 0.17602622 \\
		
		90 &  0.18192180 &  0.18192180 &  & 900 & 0.17602593 & 0.17602593 \\
	
	\hline
	\end{tabular}
	
	\end{center}
	\caption{Distribution function $F_5(k)$}
	\label{Distribution function F(k)_5(k)}
	\end{table}   

By putting zeros for $n(i,j)$ when $b > a$, we produce exactly the results that Terras~\cite{terras} obtained for the distribution function $F(k)$. Recall that we met the condition $b > a$ ($\lambda < 1$) during the analysis of the equations giving the integers $y$ (end of trajectories) as a function of $x$ (start of trajectories) in the theorem~\ref{Distribution_x_vs_y}. We deduced that there could be integers such as $x < y$ when $b > a$. As this number of integers is finite and the calculation of the distribution is carried out on all positive integers, the distribution of integers $x < y$ for $b > a$ (if it not zero) becomes negligible compared to the one where $x \ge y$. Likewise, Terras argues that the distribution of integers with $x < y$ for $b > a$ is small compared to the one where $x \ge y$ or, maybe zero. Lagarias~\cite{lagarias} makes a similar analysis. 

\section{Property of the distribution function $F(k)$}

The inherent properties of the distribution function flow directly from the properties the binomial distribution of integers (Pascal triangles) and the diophantine equations linking them (via the three theorems). We recall the fact that the column number $j$ corresponds to the exponent $k$ of the parameter $b = 2^k$, so the number of iterations, and the row number $i$ to the exponent $k_2$ of the parameter $a = 3^{k_2}$, where $k_2$ is the number of transformations on odd integers. We have $k_2 = 0, 1, \dots, k$.

If $n(i,j) = 0$ for a given combination $i = k_2$ and $j = k$, then $n(i,j) = 0$ for $i = k_2$ fixed and $j > k$. For example, $n(0,1) = 0$ ($b > a$) implies that $n(0, 2)$, $n(0, 3)$, $\cdots$, equal to $0$, because for each new value of $k$ the parameter $b$ is the previous one multiplied by $2$. The parameter $b = 2^k$ increases while the parameter $a = 3^{k_2} = 3^0$ remains constant, implying that $b$ is always greater than $a$ and $x \ge y$ (by the theorem~\ref{Distribution_x_vs_y} and the property specified at the end of the previous section). Then $\chi < k$ and the new $n(i,j) = 0$.

If $n(i,j) \ne 0$ for a given combination $i = k_2$ and $j = k$, then $n(i,j) \ne 0$ for $i > k_2$ and $j = k$ fixed. For example, $n(6,9) = 12 \ne 0$ ($b < a$) implies that $n(7, 9)$, $n(8, 9)$ and $n(9, 9)$ are different from $0$, because for each new value of $k_2$ the parameter $a$ is the previous one multiplied by $3$. The parameter $a = 3^{k_2} = 3^6$ increases while the parameter $b = 2^k = 2^9$ remains constant, implying that $b$ is always smaller than $a$ and $x < y$ (by the theorem~\ref{Distribution_x_vs_y}). Then $\chi > k$ and the new $n(i,j) \ne 0$.

Now let's look at the possible cases generated by the following two conditions, so $n(i,j) = 0$ and $n(i + 1,j) \ne 0$. Then, $n(i + 1,j + 1) = 0$ or $n(i + 1,j + 1) \ne 0$. The table~\ref{PascalTriangle_nbr_partial} represents examples of these 2 cases. 

\begin{table}[H]
\begin{center}
\begin{tabular}{c|ccccccc}
	
	$i \setminus j$ &  & 4 & 5  & 6  & 7  &  & $\cdots$ \\
		\hline

  	   &    &    &     &    &    &    &      \\
	2	   &    & 0 \cellcolor[gray]{0.6} &  0  & 0  & 0  &    & $\cdots$ \\
	3	   &    & 2 \cellcolor[gray]{0.6} &  0 \cellcolor[gray]{0.8}  & 0  & 0  &    & $\cdots$ \\
	4	   &    & 1  &  3 \cellcolor[gray]{0.8} & 3  & 0  &    & $\cdots$ \\
	5	   &    &    &  1  & 4  & 7  &    & $\cdots$ \\
		   &    &    &     &    &    &    &     \\
		   &    &    &     &    &    &    & $\cdots$    \\	

\end{tabular}
\end{center}
\caption{Pascal triangle - Partial view of the table~\ref{PascalTriangle_nbr} with $\chi$ greater than the number of iterations $k$}
\label{PascalTriangle_nbr_partial}

\end{table}

Indeed, the first condition ($n(i,j) = 0$) implies that $b/a = 2^j/3^i > 1$ and the second condition ($n(i + 1,j) \ne 0$) implies that $b/a = 2^j/3^{i+1} < 1$. By combining these two conditions, we write

\begin{equation*}
	1 < \frac{2^j}{3^i} <3.	
\label{condition_1}
\end{equation*}

Then, the quotient $b/a = 2^{j+1}/3^{i+1}$ for $n(i + 1,j + 1)$ must meet the condition 

\begin{equation*}
	\frac{2}{3} < b/a = \frac{2^j}{3^i} \cdot \frac{2}{3} < 2,	
\label{condition_2}
\end{equation*}

leading to two cases, so $n(i + 1,j + 1) = 0$ if $b/a > 1$ or $n(i + 1,j + 1) \ne 0$ if $b/a < 1$.

The first case represented by the example $n(i=2, j=4) = 0$ and $n(i+1=3, j=4) = 2$ leads to $n(i+1=3, j+1=5) = 0$.


The quotient $b/a = 2^{j+1}/3^{i+2}$ for $n(i + 2,j + 1)$ must meet the condition 

\begin{equation*}
	\frac{2}{3} \cdot \frac{1}{3} < \frac{2^j}{3^i} \cdot \frac{2}{3} \cdot \frac{1}{3} < \frac{2}{3}.	
\label{condition_3}
\end{equation*}

Then, $b < a$ and $n(i + 2,j + 1) \ne 0$, so $n(4,5) = 3$.

The second case represented by the example $n(i=3, j=5) = 0$ and $n(i+1=4, j=5) = 3$ leads to $n(i+1=4, j+1=6) = 3$.

The quotient $b/a = 2^{j+2}/3^{i+1}$ for $n(i + 1,j + 2)$ must meet the condition 

\begin{equation*}
	\frac{4}{3} < \frac{2^j}{3^i} \cdot \frac{2}{3} \cdot \frac{2}{1} < 4.	
\label{condition_4}
\end{equation*}

Then, $b > a$ and $n(i + 1,j + 2) = 0$, so $n(4, 7)= 0$.

These properties will allow us to follow the evolution of the distribution function $F_{new}(k)$ and de facto, the total number of integers $n$ with $\chi$ greater than the number of iterations $k$.   

We have $F(k) \le F(k-1)$, where we simplify the notation by using $F$ instead $F_{new}$, without losing sight of the fact the distribution function applies to integers satisfying the condition that $\chi > k$. 

Indeed, if the first non-zero value of $n(i,k-1)$ for a given $(k-1)$ is $a$, the second $b$, the third $c$, $\cdots$, and using the fact that $n(i,k) = n(i-1,k-1) + n(i,k-1)$, we have 

\begin{equation}
	F(k-1) = \frac{(a+b+c+\cdots)}{2^{k-1}}.	
\label{dist_1}
\end{equation}

In the first case,

\begin{equation}
	F(k) = \frac{(a+b)+(b+c)+\cdots)}{2^k} = \frac{2(a+b+c+\cdots)-a}{2 \cdot 2^{k-1}} = F(k-1) - \frac{a}{2^k} < F(k-1)
\label{dist_2}
\end{equation}

and in the second case,

\begin{equation}
	F(k) = \frac{((a)+(a+b)+(b+c)+\cdots)}{2^k} = \frac{2(a+b+c+\cdots)}{2 \cdot 2^{k-1}} = F(k-1).
\label{dist_3}
\end{equation}

The distribution of positive integers $F(k)$, which can simply be called density, decreases constantly without, however, reaching the zero value.

Another interesting property in Pascal's triangle (table~\ref{PascalTriangle_nbr}) is the one related to $n(k,k)$, and containing only the transformations on the odd integers. The first integer (the smallest) of a trajectory is $(2^k-1)$ and the last is $(3^k-1)$. This trajectory is 

\begin{equation*}
	(2^k-1, \phantom{1} 3 \cdot 2^{k-1}-1, \phantom{1} 3^2 \cdot 2^{k-2}-1, \phantom{1} 3^3 \cdot 2^{k-3}-1, \cdots, \phantom{1} 3^k-1).
\label{sequence_odds}
\end{equation*}

According to the theorem~\ref{periodicity} (periodicity), all integers $2^k - 1 + q \cdot 2^k$ start identical trajectories (containing only the transformations on the odd integers) of length $k + 1$.  

The distribution function for a given number of iterations (equation $\ref{dist_1}$) for $k - 1$ generates the two possible distributions (equations $\ref{dist_2}$ and $\ref{dist_3}$) for the next $k$. These results highlight two situations that will have important implication for the conclusion we will draw from them. 

A first observation is that the distribution function $F(k)$ constantly decreases with the number of iterations $k$. Terras~\cite{terras} concludes that $F(k)$ go towards zero as $k$ tends to infinity and Lagarias~\cite{lagarias} proves that the decrease is done exponentially. 

Note $N_{\chi > k}$ the number of integers and $n_{\chi > k}$ each of the integers included in the interval of $2^k$ consecutive integers that start trajectories satisfying the condition that the stopping time $\chi$ is greater than the number of iterations $k$. According to the theorem~\ref{periodicity}, the number $N_{\chi > k}$ is the same for any sequence of $2^k$ consecutive integers. Then, we add to integers $n_{\chi > k}$ all those obtained from $n_{\chi > k} + q \cdot 2^k$, where $q$ is any integer, positive, negative or zero. The numerator of the equation $\ref{dist_1}$ represent $N_{\chi > k - 1}$ for a given number $k - 1$ of iterations.  The numerators of the equations $\ref{dist_2}$ and  $\ref{dist_3}$ give $N_{\chi > k}$ after $k$ iterations. $N_{\chi > k}$ is doubled (second case analyzed previously,  equation~\ref{dist_3}) or a little less (first case, equation~\ref{dist_2}) relatively to $N_{\chi > k - 1}$. Therefore, the number of integers satisfying $\chi>k$ continuously increases with the number of iterations $k$. As $k$ tends to infinity, $N_{\chi > k}$ tends to infinity. Furthermore, $2^k$ also tends to infinity. The table~\ref{Distribution function F(k)_3(k)_with_numerator} gives the first values of $N_{\chi > k}$ and $2^k$ as $k$ is increasing. 

At first glance, the fact that $N_{\chi > k}$ tends to infinity seems to contradict the fact that the distribution function $F(k)$ tends to zero. We will see that it is not.

Indeed, $F(k)$ is given by the equation~\ref{Function_distribution_new_sum} and corresponds to the ratio $\infty / \infty$. This is similar to the quotients of continuous functions for which we seek to solve an indeterminacy by the rules of the Hospital. There are cases where the limit is finited and equal to zero. There is therefore nothing which prevents the ratio $N_{\chi > k} / 2^k = \infty / \infty$ from tending towards zero as $k$ tends towards infinity.

Because 


\begin{equation}
	\lim\limits_{k \to \infty}F(k)  \to 0  \phantom{12} or \phantom{12} \lim\limits_{k \to \infty}G(k)  \to 1
\label{limit_1}
\end{equation}

Lagarias~\cite{lagarias} concludes that almost all integers have a finite stopping time. 

We have proved that the number of integers such that $\chi > k$ has the following behavior

\begin{equation}
	\lim\limits_{k \to \infty}N_{\chi>k}(k)  \to \infty.
\label{limit_}
\end{equation}     


With this new light we can no longer conclude that almost all integers have a finite stopping time, on the contrary. Therefore, regardless of the number of iterations $k$ selected, there will always be integers that start trajectories such as $\chi > k$ in the interval of $2^k$ consecutive integers, and their number is constantly increasing with $k$. Moreover, since the density $F(k)$ of the trajectories decreases with $k$, the smallest integer starting trajectories with $\chi > k$ is possibly more and more large.  

We will conclude this paper by briefly discussing two specific points.

First we could have dealt with the possibility that there are non-trivial cycles. We have seen in the theorem~\ref{Distribution_x_vs_y} that there are possible cycles when the condition $b > a$ is satisfied. Moreover, if we examine the sequence extending from $1$ to $2^k$ positive consecutive integers instead of the one starting with $2$, we include the trivial cycle starting with the integer $1$ whose trajectory is $<1,2,1>$. The diophantine equation is written as $4y - 3x = 1$ with $b=4 > a=3$, and $y_0 = x_0 = 1$ a particular solution. The general solution is $x =  1 + 4q$ and $y =  1 + 3q$. In the table~\ref{PascalTriangle_nbr} this corresponds to $n(k_2=1,k=2)$. We can do the same for the other trivial cycles, so starting with $0$, $-1$, $-5$ and $-17$. If there are other cycles for positive integers, these must also satisfy the condition $b > a$. Anyway, if this were the case, it goes without saying that the $3x+1$ conjecture would not longer hold.

Finally, it is relatively easy to analyze the data of $F_5(k)$ generated by the iterative function $T_5$ instead of $T_3$, knowing that the three theorems used in this paper are also valid for this function.

\begin{table}
\begin{center}
\small

\rotatebox{90}{
\begin{tabular}{|c|c|c|c|c|c|c}
	\hline
	 $k$ &   number of integers $N_{\chi > k}$ & $2^k$ & $F_3(k)$  \\
	\hline

$\cdots$ &    &  &   \\

		10 &  64                   & $2^{10} = 1 \phantom{1} 024$ & $6.25 \times 10^{-2}$  \\

		20 &  $27 \phantom{1} 328$ & $2^{20} = 1 \phantom{1} 048 \phantom{1} 576$ & $2.6062 \times 10^{-2}$ \\	

	  30 &  $12 \phantom{1} 771 \phantom{1} 274$ & $2^{30} = 1 \phantom{1} 073 \phantom{1} 741 \phantom{1} 824$ & $1.1894 \times 10^{-2}$  \\	
		
		40 &  $6 \phantom{1} 402 \phantom{1} 835 \phantom{1} 000$ & $2^{40} = 1 \phantom{1} 099 \phantom{1} 511 \phantom{1} 627 \phantom{1}  776$ &  $5.8233 \times 10^{-3}$  \\
			
		50 &  $3 \phantom{1} 734 \phantom{1} 259 \phantom{1} 929 \phantom{1} 440$ & $2^{50} = 1 \phantom{1} 125 \phantom{1} 899 \phantom{1}  906 \phantom{1} 842 \phantom{1} 624$ & $3.3167 \times 10^{-3}$  \\
		
		60 &  $2 \phantom{1} 216 \phantom{1} 134 \phantom{1} 944 \phantom{1} 775 \phantom{1} 156$ & $2^{60} = 1 \phantom{1} 152 \phantom{1} 921 \phantom{1}  504 \phantom{1} 606 \phantom{1} 846 \phantom{1} 976$ &  $1.9222 \times 10^{-3}$  \\
		
		70 &  $1 \phantom{1} 241 \phantom{1} 503 \phantom{1} 538 \phantom{1} 986 \phantom{1} 719 \phantom{1} 152$ & $2^{70} = 1 \phantom{1} 180 \phantom{1} 591 \phantom{1}  620 \phantom{1} 717 \phantom{1} 411 \phantom{1} 303 \phantom{1} 424$ &  $1.0516 \times 10^{-3}$  \\
		
		80 &  $803 \phantom{1} 209 \phantom{1} 913 \phantom{1} 882 \phantom{1} 910 \phantom{1} 595 \phantom{1} 105$ & $2^{80} = 1 \phantom{1} 208 \phantom{1} 925 \phantom{1}  819 \phantom{1} 614 \phantom{1} 629 \phantom{1} 174 \phantom{1} 706 \phantom{1} 176$ & $6.6440 \times 10^{-4}$  \\
		
		90 &  $508 \phantom{1} 520 \phantom{1} 069 \phantom{1} 189 \phantom{1} 622 \phantom{1} 659 \phantom{1} 715 \phantom{1} 764$ & $2^{90} = 1 \phantom{1} 237 \phantom{1} 940 \phantom{1}  039 \phantom{1} 285 \phantom{1} 380 \phantom{1} 274 \phantom{1} 899 \phantom{1} 124 \phantom{1} 224$ & $4.1078 \times 10^{-4}$  \\
		
		100 & $302 \phantom{1} 560 \phantom{1} 669 \phantom{1} 500 \phantom{1} 543 \phantom{1} 257 \phantom{1} 546 \phantom{1} 172 \phantom{1} 187$ & $2^{100} = 1 \phantom{1} 267 \phantom{1} 650 \phantom{1} 600 \phantom{1} 228 \phantom{1} 229 \phantom{1} 401 \phantom{1} 496 \phantom{1} 703 \phantom{1} 205 \phantom{1}376$ & $2.3868 \times 10^{-4}$ \\

$\cdots$ &    &  &   \\

200 & $4.9179 \times 10^{54}$ & $2^{200} = 1.6069 \times 10^{60}$ & $3.0604\times 10^{-6}$ \\

300 & $1.1136 \times 10^{83}$ & $2^{300} = 2.0370 \times 10^{90}$ & $5.4667\times 10^{-8}$ \\

400 & $2.9920 \times 10^{111}$ & $2^{400} = 2.5822 \times 10^{120}$ & $1.1587\times 10^{-9}$ \\

500 & $8.7021 \times 10^{139}$ & $2^{500} = 3.2734 \times 10^{150}$ & $2.6584\times 10^{-11}$ \\

600 & $2.6746 \times 10^{168}$ & $2^{600} = 4.1495 \times 10^{180}$ & $6.4455\times 10^{-13}$ \\

700 & $8.2683 \times 10^{196}$ & $2^{700} = 5.2601 \times 10^{210}$ & $1.5719\times 10^{-14}$ \\

800 & $2.7314 \times 10^{225}$ & $2^{800} = 6.6680 \times 10^{241}$ & $4.0963\times 10^{-16}$ \\

900 & $9.1605 \times 10^{253}$ & $2^{900} = 8.4527 \times 10^{270}$ & $1.0837\times 10^{-17}$ \\

$\cdots$ &    &  &   \\
	
	\hline
	\end{tabular}
	}
	\end{center}
	\caption{Number of integers starting trajectories with $\chi>k$ in the interval of $2^k$ consecutive integers, different total number of trajectories ($2^k$) for $k$ iterations, and distribution function $F_3(k)$}
	\label{Distribution function F(k)_3(k)_with_numerator}
	\end{table}

\end{document}